# Large sample asymptotics for the
# two-parameter Poisson–Dirichlet process

## Lancelot F. James[*1]

*Hong Kong University of Science and Technology*

**Abstract:** This paper explores large sample properties of the two-parameter
$(\alpha, \theta)$ Poisson–Dirichlet Process in two contexts. In a Bayesian context of es-
timating an unknown probability measure, viewing this process as a natural
extension of the Dirichlet process, we explore the consistency and weak conver-
gence of the the two-parameter Poisson–Dirichlet posterior process. We also
establish the weak convergence of properly centered two-parameter Poisson–
Dirichlet processes for large $\theta + n\alpha$. This latter result complements large $\theta$
results for the Dirichlet process and Poisson–Dirichlet sequences, and comple-
ments a recent result on large deviation principles for the two-parameter
Poisson–Dirichlet process. A crucial component of our results is the use of
distributional identities that may be useful in other contexts.

## Contents

## 1. Introduction

In this work, for $0 \le \alpha < 1$ and $\theta > -\alpha$, we are interested in the two-parameter
class of random probability measures that are formed by

$$P_{\alpha,\theta}(\cdot) \stackrel{d}{=} \sum_{k=1}^{\infty} V_k \prod_{j=1}^{k-1} (1 - V_j) \delta_{Z_k}(\cdot) \qquad (1.1)$$

where the $V_k$ are independent beta$(1-\alpha, \theta+k\alpha)$ random variables and, independent
of these, the $(Z_k)$ are an i.i.d. sequence with values in some Polish space $I$ with com-
mon (non-atomic) distribution $H$. That is to say $P_{\alpha,\theta}$ is a random probability mea-
sure taking values in $\mathcal{P}_I$, where $\mathcal{P}_I$ is the set of all probability measures on $I$. We will

---

*Supported in part by the Grants RGC-HKUST 600907, SBM06/07.BM14 and
HIA05/06.BM03 of the HKSAR.

[1]Lancelot F. James, Hong Kong University of Science and Technology, Department of Informa-
tion and Systems Management, Clear Water Bay, Kowloon, Hong Kong, e-mail: lancelot@ust.hk.

*AMS 2000 subject classifications:* Primary 62G05; secondary 62F15.

*Keywords and phrases:* Bayesian consistency, multiplier CLT, two-parameter Poisson–Dirichlet
process, weak convergence.

simply say that a random probability measure $P$ is a two-parameter $(\alpha, \theta)$ Poisson–Dirichlet process, having law say $\Pi_{\alpha,\theta}$ on $\mathcal{P}_I$, suppressing dependence on $H$, if $P$ can be represented as in (1.1). That is $P(\cdot) \overset{d}{=} P_{\alpha,\theta}(\cdot)$. For shorthand we write $P \sim \Pi_{\alpha,\theta}$. We denote the expectation operator corresponding to the law $\Pi_{\alpha,\theta}$ as $E_{\alpha,\theta}$, which is such that $E_{\alpha,\theta}[P(\cdot)] = H(\cdot)$. We note that the $V_k$ are obtained by size-biasing a ranked sequence of probabilities known as the two-parameter Poisson–Dirichlet sequence. It follows that by permutation invariance $P_{\alpha,\theta}$ may also be represented in terms of this sequence. Many properties of the two-parameter Poisson–Dirichlet sequence, as it primarily relates to Bessel and Brownian phenomena, were discussed in [36]. This sequence has gained in importance as it is seen to arise in a number of different areas including, for instance, Bayesian statistics, population genetics and random fragmentation and coalescent theory connected to physics. See [34] for some updated references and [19] for some connections to Dirichlet means.

When $\alpha = 0$ then $P$ is a Dirichlet process in the sense of Ferguson [9]. Pitman [35] showed that within a Bayesian context, these random probability measures can be seen as natural and quite tractable extensions of the Dirichlet process. In particular Pitman [35] showed that if random variables $X_1, \ldots, X_n$ given $P$ are i.i.d. $P$ and $P$ has prior distribution $\Pi_{\alpha,\theta}$, then the posterior distribution of $P|X_1, \ldots, X_n$, denoted as $\Pi_{\alpha,\theta}^{(n)}$, corresponds to the law of the random probability measure,

$$P_{\alpha,\theta}^{(n)}(\cdot) = R_{n(\mathbf{p})} P_{\alpha,\theta+n(\mathbf{p})\alpha}(\cdot) + (1 - R_{n(\mathbf{p})}) D_n(\cdot)$$

where

$$D_n(\cdot) = \sum_{j=1}^{n(\mathbf{p})} \Delta_j \delta_{Y_j}(\cdot).$$

$(\Delta_1, \ldots, \Delta_{n(\mathbf{p})})$ is a Dirichlet$(e_1 - \alpha, \ldots, e_{n(\mathbf{p})} - \alpha)$ random vector. All the random variables appearing on the right hand side are conditionally independent given the data. $R_{n(\mathbf{p})}$ is a beta$(\theta + n(\mathbf{p})\alpha, n - n(\mathbf{p})\alpha)$ random variable and $\{Y_1, \ldots, Y_{n(\mathbf{p})}\}$ denotes the $1 \leq n(\mathbf{p}) \leq n$ unique values of $\{X_1, \ldots, X_n.\}$ Furthermore, $P_{\alpha,\theta+n(\mathbf{p})\alpha} \sim \Pi_{\alpha,\theta+n(\mathbf{p})\alpha}$. $e_j$ is the number of $X_i$ equivalent to $Y_j$ for $j = 1, \ldots, n(\mathbf{p})$., When $\alpha = 0$ one obtains the posterior distribution derived in Ferguson [9]. The notation $e_j$ and $n(\mathbf{p})$ are taken from Lo [29], as discussed in Ishwaran and James [17]. Furthermore, a generalization of the Blackwell–MacQueen [4] prediction rule is given by

$$(1.2) \qquad P(X_{n+1} \in \cdot \,| X_1, \ldots, X_n) = \frac{\theta + n(\mathbf{p})\alpha}{\theta + n} H(\cdot) + \sum_{j=1}^{n(\mathbf{p})} \frac{(e_j - \alpha)}{\theta + n} \delta_{Y_j}(\cdot).$$

Note also that

$$P(X_{n+1} \in \cdot \,| X_1, \ldots, X_n) = E[P_{\alpha,\theta}^{(n)}(\cdot)].$$

We also write

$$\tilde{F}_n(\cdot) = \sum_{j=1}^{n(\mathbf{p})} \frac{(e_j - \alpha)}{n - n(\mathbf{p})\alpha} \delta_{Y_j}(\cdot)$$

which importantly reduces to the empirical distribution when $\alpha = 0$, or when $n(\mathbf{p}) = n$.

Let $P_0^\infty$ denote a product measure on $I^\infty$ making $X_i$ for $i = 1, \ldots, \infty$ independent with common (true) distribution $P_0$. In this paper, extending the known case of the Dirichlet process, we show that as $n \to \infty$ the posterior distribution $\Pi_{\alpha,\theta}^{(n)}$ behaves as follows. When $P_0$ is discrete then $\Pi_{\alpha,\theta}^{(n)}$ converges weakly to a point mass at $P_0$ a.s. $P_0^{\infty}$. When $P_0$ is continuous, that is, non-atomic, $\Pi_{\alpha,\theta}^{(n)}$ converges weakly to a point mass at the mixture $\alpha H + (1 - \alpha)P_0$. Thus when $P_0$ is discrete/atomic the posterior distribution is consistent. However, when $P_0$ is non-atomic the posterior distribution is inconsistent, unless either $\alpha = 0$ which corresponds to the case of the Dirichlet process, or more implausibly one chooses $H = P_0$. In addition to this result we establish a functional central limit theorem, for the case where $P_0$ is non-atomic, by showing that the process $P_{\alpha,\theta}^{(n)}$ centered at its expectation (1.2), indexed over classes of functions, converges weakly to a Gaussian process. This is in line with nonparametric Bernstein–Von Mises results of for instance, [31], [30], [28], [5] and [22]. Additionally, we note that the weak convergence of the two-parameter Poisson–Dirichlet process may be of interest in other fields. In particular, in order for us to discuss the posterior weak convergence, when $P_0$ is non-atomic, we will need to address the weak convergence of the centered process

$$\nu_{\alpha,\theta+n\alpha}(\cdot) = \sqrt{n}(P_{\alpha,\theta+n\alpha} - H)(\cdot),$$

as $n \to \infty$, which poses additional challenges. Note this process does not depend on the data, except through the sample size $n$. Furthermore, the study of $\nu_{\alpha,\theta+n\alpha}$ is more in line with the literature on the behavior of Dirichlet processes and Poisson–Dirichlet sequences when $\theta \to \infty$. See, for instance, [21], [7] and [32]. Additionally, our work is complementary to a recent result of [8] on large deviation principles for the two-parameter Poisson–Dirichlet process.

Returning to the consistency result, in terms of estimating the true $P_0$ in a nonparametric statistics setting, our result shows that unlike the case of the Dirichlet process, it is perhaps unwise to use $\Pi_{\alpha,\theta}$ as a prior. However we should point out that Ishwaran and James [17, 18], owing to the attractive results in [35], suggested that one could use $\Pi_{\alpha,\theta}$ in a mixture modeling setting analogous to the case of the Dirichlet process in Lo [29]. In this more formidable setting one can deduce strong consistency of the posterior density, induced by the priors $\Pi_{\alpha,\theta}$ on the mixing distribution, by using the results of Lijoi, Prünster and Walker [25]. In fact, in a work subsequent to ours, this was recently shown by Jang, Lee and Lee [20]. We also note that those authors also obtain our consistency result as a special case. For some more results on the modern treatment of Bayesian consistency in nonparametric settings one may note, for instance, the works of Ghosal, Ghosh and Ramamoorthi [11], Barron, Schervish and Wasserman [2] and Ghosal, Ghosh and van der Vaart [12]; and the book of Ghosh and Ramamoorthi [13].

## 2. Consistency

This section describes the consistency behavior of the posterior distribution in the case where the true distribution $P_0$ is either continuous or discrete. First we note the following fact;

**Lemma 2.1.** *Let $f$ and $g$ denote measurable functions on $I$, then for $0 \le \alpha < 1$ and $\theta > -\alpha$,*

$$E_{\alpha,\theta}[P(f)P(g)] = \frac{\theta + \alpha}{\theta + 1}H(g)H(f) + \frac{1 - \alpha}{\theta + 1}H(fg).$$

*Proof.* The proof proceeds by using disintegrations. First the joint distribution of $(X_1, P)$ can be written as,

$$P(dx_1)\Pi_{\alpha,\theta}(dP) = \Pi_{\alpha,\theta}(dP|x_1)H(dx_1)$$

where $\Pi_{\alpha,\theta}(dP|x_1) = \Pi_{\alpha,\theta}^{(1)}(dP)$. Then,

$$P(dx_2)\Pi_{\alpha,\theta}(dP|x_1) = \Pi_{\alpha,\theta}^{(2)}(dP)E[P(dx_2)|X_1 = x_1]$$

where $\Pi_{\alpha,\theta}^{(2)}(dP) = \Pi_{\alpha,\theta}(dP|x_1, x_2)$ and

$$E[P(dx_2)|X_1 = x_1] = \frac{\theta+\alpha}{\theta+1}H(dx_2) + \frac{1-\alpha}{\theta+1}\delta_{x_1}(dx_2).$$

Now this gives

$$E_{\alpha,\theta}[P(dx_1)P(dx_2)] = H(dx_1)\left[\frac{\theta+\alpha}{\theta+1}H(dx_2) + \frac{1-\alpha}{\theta+1}\delta_{x_1}(dx_2)\right].$$

which by writing $P(g)P(f) = \int_I \int_I g(x_1)f(x_2)P(dx_1)P(dx_2)$ completes the result. □

We proceed as in Diaconis and Freedman[10] by showing that the posterior distribution concentrates around the prediction rule. First using Diaconis and Freedman ([10], p. 1117), we define a suitable class of semi-norms such that convergence under such norms implies weak convergence. Let $\mathcal{A} = \bigcup_{i=1}^{\infty} A_i$ be a partition of $I$. Then define the semi-norm between probability measures

$$(2.1) \qquad \mid P - Q \mid_{\mathcal{A}} = \sqrt{\sum_{i=1}^{\infty}\left[P(A_i) - Q(A_i)\right]^2},$$

for a suitable generating sequence of partitions $\mathcal{A}$ where, naturally for any $Q$, $\sum_{i=1}^{\infty} Q(A_i) = 1$. Now similar to ([10], Equation 14) for the Dirichlet process, we will show that the posterior distribution concentrates around the prediction rule (1.2).

In order to do this one only needs to evaluate the posterior expectation, expressible as,

$$(2.2) \qquad E\left[\mid P_{\alpha,\theta}^{(n)} - E[P_{\alpha,\theta}^{(n)}]\mid_{\mathcal{A}}^2\right]$$

where $E[P_{\alpha,\theta}^{(n)}]$ equates with the prediction rule probability given in (1.2). We obtain,

**Lemma 2.2.**
$$E\left[\mid P_{\alpha,\theta}^{(n)} - E[P_{\alpha,\theta}^{(n)}]\mid_{\mathcal{A}}^2\right] \leq \frac{1}{\theta+n+1}.$$

*Proof.* First using basic ideas we expand, for each set $A_i$,

$$\left(P_{\alpha,\theta}^{(n)}(A_i) - E[P_{\alpha,\theta}^{(n)}(A_i)]\right)^2.$$

Furthermore,

$$\begin{aligned}(P_{\alpha,\theta}^{(n)}(A_i))^2 &= R_{n(\mathbf{p})}^2 P_{\alpha,\theta+n(\mathbf{p})\alpha}^2(A_i) \\ &+ 2R_{n(\mathbf{p})}(1 - R_{n(\mathbf{p})})P_{\alpha,\theta+n(\mathbf{p})\alpha}(A_i)D_n(A_i) + (1 - R_{n(\mathbf{p})})^2 D_n^2(A_i).\end{aligned}$$

Now from Lemma 2.1

$$E[P^2_{\alpha, \theta + n(\mathbf{p})\alpha}(A_i)] = \frac{\theta + n(\mathbf{p})\alpha + \alpha}{\theta + n(\mathbf{p})\alpha + 1} H^2(A_i) + \frac{1 - \alpha}{\theta + n(\mathbf{p})\alpha + 1} H(A_i).$$

Additionally,

$$E[R^2_{n(\mathbf{p})}] = \frac{(\theta + n(\mathbf{p})\alpha)(\theta + n(\mathbf{p})\alpha + 1)}{(\theta + n)(\theta + n + 1)},$$

$$E[(1 - R_{n(\mathbf{p})})^2] = \frac{(n - n(\mathbf{p})\alpha)(n - n(\mathbf{p})\alpha + 1)}{(\theta + n)(\theta + n + 1)}$$

and

$$E[\Delta_l \Delta_j] = \frac{(e_j - \alpha)(e_l - \alpha)}{(n - n(\mathbf{p})\alpha)(n - n(\mathbf{p})\alpha + 1)}.$$

It follows that

$$\begin{aligned}
E[(P^{(n)}_{\alpha, \theta}(A_i))^2] &= \frac{(\theta + n(\mathbf{p})\alpha)(\theta + n(\mathbf{p})\alpha + \alpha)}{(\theta + n)(\theta + n + 1)} H^2(A_i) \\
&+ 2\frac{(\theta + n(\mathbf{p})\alpha)(n - n(\mathbf{p})\alpha)}{(\theta + n)(\theta + n + 1)} \tilde{F}_n(A_i) H(A_i) \\
&+ \frac{(n - n(\mathbf{p})\alpha)^2}{(\theta + n)(\theta + n + 1)} \tilde{F}^2_n(A_i).
\end{aligned}$$

Now

$$\begin{aligned}
(E[P^{(n)}_{\alpha, \theta}(A_i)])^2 &= \frac{(\theta + n(\mathbf{p})\alpha)^2}{(\theta + n)^2} H^2(A_i) \\
&+ 2\frac{(\theta + n(\mathbf{p})\alpha)(n - n(\mathbf{p})\alpha)}{(\theta + n)^2} \tilde{F}_n(A_i) H(A_i) \\
&+ \frac{(n - n(\mathbf{p})\alpha)^2}{(\theta + n)^2} \tilde{F}^2_n(A_i).
\end{aligned}$$

Taking differences and using the fact that $\sum_{i=1}^{\infty} F_n(A_i) H(A_i) \leq 1$ and similar arguments completes the result. □

**Proposition 2.1.** *If $P_0$ is continuous then the posterior distribution $\Pi^{(n)}_{\alpha, \theta}$ converges weakly to point mass at the distribution*

$$\alpha H(\cdot) + (1 - \alpha) P_0(\cdot) \ a.e. \ P^{\infty}_0.$$

*Hence the posterior is consistent only if either $P$ is a Dirichlet process or $H = P_0$.*

*Proof.* From Lemma 2.2 it follows that the posterior distribution must concentrate around the prediction rule. Now, under $P_0$ (assuming a continuous $P_0$) the prediction rule becomes,

$$(2.3) \qquad P(X_{n+1} \in \cdot \mid X_1, \ldots, X_n) = \frac{\theta + n\alpha}{\theta + n} H(\cdot) + \sum_{j=1}^{n} \frac{(1 - \alpha)}{\theta + n} \delta_{X_j}(\cdot).$$

It is clear that under $P_0$, that $P(X_{n+1} \in \cdot \mid X_1, \ldots, X_n)$ converges uniformly over appropriate Glivenko–Cantelli classes to

$$\alpha H(\cdot) + (1 - \alpha) P_0(\cdot)$$

for almost all sample sequences. One gets this by simple algebra and classical results about empirical processes (the second term in (2.3)) appropriately modified. □

The previous results says that the posterior distribution is inconsistent for all non-atomic $P_0$ unless $\alpha = 0$ or one has chosen $H = P_0$. The behavior in the case where $P_0$ admits ties is quite different and is summarized in the next result

**Proposition 2.2.** *Suppose that $P_0$ is a discrete law such that $n(\mathbf{p})/n \to 0$ then the posterior distribution $\Pi_{\alpha,\theta}^{(n)}$ converges weakly to point mass at $P_0$, a.e. $P_0^\infty$, for all $0 \le \alpha < 1$ and $\theta > -\alpha$.*

*Proof.* Since, $n(\mathbf{p})/n \to 0$, it follows that

$$(2.4) \qquad P(X_{n+1} \in \cdot \mid X_1, \ldots, X_n) = \frac{\theta + n(\mathbf{p})\alpha}{b+n} H(\cdot) + \sum_{j=1}^{n(\mathbf{p})} \frac{(e_j - \alpha)}{\theta + n} \delta_{Y_j}(\cdot),$$

converges uniformly to $P_0$ for almost all sample sequences $X_1, X_2, \ldots,$. This is true since,

$$\sum_{j=1}^{n(\mathbf{p})} \frac{(e_j - \alpha)}{\theta + n} \delta_{Y_j} = \sum_{i=1}^{n} \frac{1}{\theta + n} \delta_{X_i} - \sum_{j=1}^{n(\mathbf{p})} \frac{\alpha}{\theta + n} \delta_{Y_j}$$

where the second term on the right converges to zero. $\qquad \square$

### *2.1. Some more limits*

Note that one obtains some information on the limit behavior of the random probability measure $P_{\alpha,\theta+n(\mathbf{p})\alpha}$ which has law $\Pi_{\alpha,\theta+n(\mathbf{p})\alpha}$. This type of result is more in line with large $\theta$ type asymptotics. In this case large $\theta$ is replaced by large $n(\mathbf{p})\alpha$.

**Proposition 2.3.** *Suppose that $P_0$ is such that $n(\mathbf{p}) \to \infty$, then the two-parameter Poisson–Dirichlet law $\Pi_{\alpha,\theta+n(\mathbf{p})\alpha}$ converges weakly to point mass at $H$.*

*Proof.* The proof proceeds by again utilizing the semi-norm in (2.1). We have that $E_{\alpha,\theta+n(\mathbf{p})\alpha}[P] = H$. Furthermore, from Lemma 2.1. the variance of $P(A)$ under $\Pi_{\alpha,\theta+n(\mathbf{p})\alpha}$ is, for each $A$, equal to,

$$\frac{1-\alpha}{\theta+n(\mathbf{p})\alpha+1} H(A)(1-H(A)).$$

Hence if $P$ has law $\Pi_{\alpha,\theta+n(\mathbf{p})\alpha}$, then,

$$E_{\alpha,\theta+n(\mathbf{p})\alpha}\left[|P-H|_{\mathcal{A}}^2\right] \le \frac{1-\alpha}{\theta+n(\mathbf{p})\alpha+1},$$

completing the result. $\qquad \square$

## 3. Bernstein–von Mises and functional central limit theorems

In this section we address the more formidable problem of establishing functional central limit theorems for centered versions of the two-parameter process. We will restrict ourselves to the case where $P_0$ is continuous which presents some difficulties. In particular, in that setting, we are interested in the asymptotic behavior in distribution, as $n \to \infty$, of the posterior process

$$\nu_{\alpha,\theta}^{(n)}(\cdot) = \sqrt{n}(P_{\alpha,\theta}^{(n)} - E[P_{\alpha,\theta}^{(n)}])(\cdot)$$

conditional on the sequence $X_1, X_2, \ldots$, and the asymptotic behavior of the process,

$$\nu_{\alpha, \theta + n\alpha}(\cdot) := \sqrt{n}(P_{\alpha, \theta + n\alpha} - H)(\cdot),$$

uniform over classes of functions $\mathcal{F}$. For clarity, we first mention some elements of the (modern) framework of weak convergence of stochastic processes on general function indexed Banach spaces. There is a rich literature on this subject, here we use as references [39], [38], Ch. 10, and [14]. Let $\mathcal{F}$ denote a collection of measurable functions $f : I \to \Re$ and let $\ell^\infty(\mathcal{F})$ denote the set of all uniformly bounded real functions on $\mathcal{F}$. Now for a random probability measure $Q_n$, and the probability measure defined as its expectation $E[Q_n] := Q$, we consider the maps from $\mathcal{F} \to \Re$ given by the linear functional

$$f \to Q_n(f) = \int_I f(x) Q_n(dx)$$

and

$$f \to Q(f) = \int_I f(x) Q(dx).$$

$\mathbb{G}_n(\cdot) = \sqrt{n}(Q_n - Q)(\cdot)$ denotes its centered process and let $\mathbb{G}_Q$ denote a Gaussian process with zero mean and covariances

$$(3.1) \qquad E[\mathbb{G}_Q(f)\mathbb{G}_Q(g)] = Q(fg) - Q(f)Q(g).$$

A Gaussian process with covariance (3.1) is said to be a $Q$-Brownian bridge. We will assume that $\mathcal{F}$ possesses enough measurability for randomization and write, as in ([37], p. 2056), $\mathcal{F} \in M(Q)$. The notation $L_2(Q)$ represents the equivalence class of square integrable functions. Furthermore, a function $F(x)$ such that $|f(x)| \leq F(x)$ for all $x$ and $f \in \mathcal{F}$ is said to be an envelope.

We are interested in the cases where the sequence of processes $\{G_n(f) : f \in \mathcal{F}\} \in \ell^\infty(\mathcal{F})$ converges in distribution to a Gaussian process $\mathbb{G}_Q$ uniformly over $\mathcal{F}$. In that case we write

$$G_n \rightsquigarrow \mathbb{G}_Q \text{ in } \ell^\infty(\mathcal{F}).$$

Furthermore, because we are interested in the convergence of posterior distributions, if $G_n$ depends on data $X_1, X_2 \ldots$, we will need the more delicate notion of conditional weak convergence of $G_n(\cdot)$ for almost all sample sequences $X_1, X_2 \ldots$, and we write

$$G_n \rightsquigarrow \mathbb{G}_Q \text{ in } \ell^\infty(\mathcal{F}) \text{ a.s. }.$$

More formally, one may say that the processes converge in the sense of a bounded Lipschitz metric outer almost surely (see [14] and [39]).

The rich theory of weak convergence of empirical processes addresses the cases where $Q_n = P_n = n^{-1} \sum_{i=1}^n \delta_{X_i}(\cdot)$ is the empirical measure, and $Q = E_{P_0}[P_n] = P_0$ is the true underlying distribution of the data. Hence one has,

$$(3.2) \qquad \sqrt{n}(P_n - P_0) \rightsquigarrow \mathbb{G}_{P_0} \text{ in } \ell^\infty(\mathcal{F})$$

for many classes of $\mathcal{F}$. The classes are said to be $P_0$-Donsker. The classical case of the empirical distribution function, $F_n(t) = \sum_{i=1}^n I_{(X_i \leq t)}/n$ arises by setting $f_t(x) = I_{(-\infty < x \leq t)}$ for $t$ ranging over $\Re$; see ([38], Example 19.6).

The most notable results for convergence conditionally on the data, center around the bootstrap and its exchangeably weighted extensions where

$$Q_n(\cdot) := P_W(\cdot) = \sum_{i=1}^\infty W_i \delta_{X_i}(\cdot)$$

for $(W_i)$ an exchangeable sequence of positive weights summing to 1. In particular, we will use the following result of ([37],Theorem 2.1),

$$(3.3) \qquad \sqrt{n}(P_W - P_n) \rightsquigarrow c\mathbb{G}_{P_0} \text{ in } \ell^\infty(\mathcal{F}), a.s.$$

provided that (3.2) holds, $P_0(F^2) < \infty$, and the weights satisfy certain conditions as given in [37]. The constant $c$ is determined by the weights. The result generalizes the result of [16] for Efron's bootstrap empirical measure. In the case of Efron's bootstrap and the Bayesian bootstrap $c = 1$. For results on the real line see [3], [27] and [33].

As mentioned at the beginning of this section, we will consider the behavior of the process $\nu_{\alpha,\theta}^{(n)}$ in the case where $P_0$ is continuous. We shall see that we can handle part of the weak convergence of $\nu_{\alpha,\theta}^{(n)}$ by utilizing results in [37]. This is very much in the spirit of [27], [26] and [5]. However we will also need to deal with the behavior of the process $\nu_{\alpha,\theta+n\alpha}$. This process is considerably more challenging to handle as it is not obviously related to an empirical-type measure. However, in Section 4, we will exploit an important distributional identity that allows us to treat $\nu_{\alpha,\theta+n\alpha}$, as a *measurelike* sum of i.i.d. processes in the sense of [1], [39], Section 2.11.1.1, and [41]. Throughout we will assume that $\mathcal{F} \in M(H)$ and that there exist an envelope $F(x)$ satisfying $H(F^2) < \infty$.

We shall also have need of the (unconditional) multiplier central limit theorem on Banach spaces for i.i.d. random variables in the Lorentz $L_{2,1}$ space, which is found in [24]. More details may be found in [39] (see also [15] and [14]). A random variable $\xi$ (see [39], Section 2.9) is said to be in $L_{2,1}$ if

$$\|\xi\|_{2,1} = \int_0^\infty \sqrt{P(|\xi| > x)}dx < \infty.$$

Finiteness of $\|\xi\|_{2,1}$ requires slightly more than a second moment but is implied by a $2 + \epsilon$ absolute moment. In our case, the $L_{2,1}$ condition is easily satisfied as the variables we shall encounter are gamma random variables having all moments.

### *3.1. Continuous case*

We now address the case of weak convergence of $\nu_{\alpha,\theta}^{(n)}$ when $P_0$ is continuous. Here we will need to consider the process $\nu_{\alpha,\theta+n\alpha}$. Because, we will not be working strictly with empirical-type processes, but actually measure-like processes, we will restrict ourselves to the quite rich class of $\mathcal{F}$ that constitute a Vapnik–Chervonenkis graph class (VCGC) (see for instance [41], p. 239). This avoids the need to otherwise state uniform entropy-type conditions. We now present the result below for $\nu_{\alpha,\theta}^{(n)}$. We will present a partial proof of this result and then address the behavior of $\nu_{\alpha,\theta+n\alpha}$, in the next section.

**Theorem 3.1.** *Let $\mathcal{F}$ be a VCGC subclass of $L_2(P_0)$ and $L_2(H)$ with envelope $F$ such that $P_0(F^2) < \infty$ and $H(F^2) < \infty$. For $0 \le \alpha < 1$, and $\theta > -\alpha$, let $\nu_{\alpha,\theta}^{(n)}(\cdot)$ denote the posterior, $\Pi_{\alpha,\theta}^{(n)}$, process centered at its mean and scaled by $\sqrt{n}$. Then when $P_0$ is continuous, conditional on the sequence $X_1, X_2, \ldots,$*

$$\nu_{\alpha,\theta}^{(n)} \rightsquigarrow \sqrt{1-\alpha}\,\mathbb{G}_{P_0} + \sqrt{\alpha(1-\alpha)}\,\mathbb{G}_H + \sqrt{\alpha}\tilde{N}(P_0 - H) \text{ in } \ell^\infty(\mathcal{F}) \text{ a.s.}$$

*Where $\mathbb{G}_{P_0}$ and $\mathbb{G}_H$ are independent Gaussian processes, independent of $\tilde{N}$ which is a standard Normal random variable.*

*Proof.* The process $\nu_{\alpha,\theta}^{(n)}$ is equivalent to

$$\sqrt{n}(1-R_n)(D_n - \tilde{F}_n) + \sqrt{n}(R_n - \frac{\theta+n\alpha}{\theta+n})(H - \tilde{F}_n) + R_n\nu_{\alpha,\theta+n\alpha}.$$

Now since $P_0$ is continuous it follows that $\tilde{F}_n$ equates with the empirical measure $P_n(\cdot) = n^{-1}\sum_{i=1}^n \delta_{X_i}(\cdot)$, and $D_n = P_W$ where $P_W$ has weights

$$W_i = \frac{G_{1-\alpha,i}}{\sum_{l=1}^n G_{1-\alpha,l}}$$

where $G_{1-\alpha,i}$ are i.i.d. gamma$(1-\alpha,1)$ random variables. Furthermore, $R_n$ is beta$(\theta+n\alpha, n(1-\alpha))$ and hence converges in probability to $\alpha$ as $n \to \infty$. So the first term is asymptotically equivalent to the process

$$\sqrt{n}(1-\alpha)(P_W - P_n)$$

and it follows that the result of ([37], Theorem 2.1) applies. Hence the process, without the $(1-\alpha)$ term, satisfies (3.3) with $c = 1/\sqrt{1-\alpha}$. It is easy to see that $R_n$ centered by its mean and scaled by $\sqrt{n}$ converges to a Normal distribution, hence the second term converges in distribution to $\sqrt{\alpha}\tilde{N}(H - P_0)$. Finally, the limit of $\nu_{\alpha,\theta+n\alpha}$ will be verified in the next section. □

## 4. Asymptotic behavior of a Poisson–Dirichlet $(\alpha, \theta + n\alpha)$ process

In this last section we establish the weak convergence of the process $\nu_{\alpha,\theta+n\alpha}$. Since $P_{\alpha,\theta+n\alpha}$ is closely associated with various properties of Brownian and Bessel processes we expect that this result will be of interest in those settings. The establishment of this result requires a few non-trivial maneuvers as $P_{\alpha,\theta+n\alpha}$ does not appear to have any similarities to an empirical process. We first establish an important distributional identity

**Proposition 4.1.** *Let $P_{\alpha,\theta+n\alpha}$ denote a random probability measure with law $\Pi_{\alpha,\theta+n\alpha}$, then*

$$P_{\alpha,\theta+n\alpha}(\cdot) \overset{d}{=} \frac{G_\theta}{G_{\theta+n\alpha}}P_{\alpha,\theta}(\cdot) + \sum_{i=1}^n \frac{G_{\alpha,i}}{G_{\theta+n\alpha}}P_{\alpha,\alpha}^{(i)}(\cdot)$$

*where $P_{\alpha,\alpha}^{(i)}$ are i.i.d. $\Pi_{\alpha,\alpha}$ random probability measures independent of $G_\theta$, $(G_{\alpha,i})$, where $G_{\alpha,i}$ are i.i.d. gamma$(\alpha,1)$ random variables, $G_\theta$ is gamma$(\theta,1)$, independent of $(G_{\alpha,i})$ and $G_{\theta+n\alpha} = G_\theta + \sum_{i=1}^n G_{\alpha,i}$.*

*Proof.* It is enough to establish this result for $P_{\alpha,\theta+n\alpha}(g)$, for an arbitrary bounded measurable function $g$. The result follows by noting that for any $\theta > 0$,

$$E[e^{-\lambda G_\theta P_{\alpha,\theta}(g)}] = \left[\int_I (1 + \lambda g(x))^\alpha H(dx)\right]^{-\theta/\alpha}$$

which is equivalent to the Cauchy–Stieltjes transform of order $\theta$ of $P_{\alpha,\theta}(g)$, whose form was obtained by [40]. It follows that,

$$E[e^{-\lambda G_{\theta+n\alpha} P_{\alpha,\theta+n\alpha}(g)}] = \left[\int_I (1 + \lambda g(x))^\alpha H(dx)\right]^{-\theta/\alpha - n}$$

which is the same as

$$E[e^{-\lambda G_\theta P_{\alpha,\theta}(g)}] \prod_{i=1}^n E[e^{-\lambda G_{\alpha,i} P_{\alpha,\alpha}^{(i)}(g)}].$$

Thus we can conclude that

$$G_{\theta+n\alpha} P_{\alpha,\theta+n\alpha} \stackrel{d}{=} G_\theta P_{\alpha,\theta} + \sum_{i=1}^n G_{\alpha,i} P_{\alpha,\alpha}^{(i)}(\cdot).$$

Now, using the fact that $G_{\theta+n\alpha} = G_\theta + \sum_{i=1}^n G_{\alpha,i}$ is gamma$(\theta + n\alpha, 1)$, it follows using the calculus of beta and gamma random variables that

$$G_{\theta+n\alpha} P_{\alpha,\theta+n\alpha}(\cdot) \stackrel{d}{=} G_{\theta+n\alpha} \left[ \frac{G_\theta}{G_{\theta+n\alpha}} P_{\alpha,\theta}(\cdot) + \sum_{i=1}^n \frac{G_{\alpha,i}}{G_{\theta+n\alpha}} P_{\alpha,\alpha}^{(i)}(\cdot) \right]$$

where on the right hand side, $G_{\theta+n\alpha}$ is independent of the term in brackets. Now we use the fact that gamma random variables are simplifiable to conclude the result. See Chaumont and Yor ([6], Sections 1.12 and 1.13) for details on simplifiable random variables. □

Proposition 4.1 now allows us to write

$$\nu_{\alpha,\theta+n\alpha} = \sqrt{n} \frac{G_\theta}{G_{\theta+n\alpha}} P_{\alpha,\theta}(\cdot) + \sqrt{n} \frac{\sum_{i=1}^n G_{\alpha,i}}{G_{\theta+n\alpha}} \sum_{i=1}^n \frac{G_{\alpha,i}}{\sum_{i=1}^n G_{\alpha,i}} P_{\alpha,\alpha}^{(i)}(\cdot) - \sqrt{n} H(\cdot).$$

Additionally one has for $\theta > 0$, the covariance formula

$$cov(\frac{G_\theta}{\theta}(P_{\alpha,\theta} - H)(f), \frac{G_\theta}{\theta}(P_{\alpha,\theta} - H)(g)) = \frac{1-\alpha}{\theta}[H(fg) - H(f)H(g)],$$

which follows from Lemma 2.1. Using these points we obtain the next result.

**Theorem 4.1.** *Let $\mathcal{F}$ be a VCGC subclass of $L_2(H)$ with envelope $F$ such that $H(F^2) < \infty$. Let, for $0 < \alpha < 1$ and $\theta > -\alpha$, $P_{\alpha,\theta+n\alpha}$ denote the random probability measure with Poisson–Dirichlet law $\Pi_{\alpha,\theta+n\alpha}$, and define $\nu_{\alpha,\theta+n\alpha}(\cdot) := \sqrt{n}(P_{\alpha,\theta+n\alpha} - H)(\cdot)$. Then as $n \to \infty$,*

$$\nu_{\alpha,\theta+n\alpha} \rightsquigarrow \frac{\sqrt{1-\alpha}}{\sqrt{\alpha}} \mathbb{G}_H \text{ in } \ell^\infty(\mathcal{F}).$$

*Proof.* The key to this result is of course Proposition 4.1, which allows us to express $\nu_{\alpha,\theta+n\alpha}$, in terms of i.i.d. processes $P_{\alpha,\alpha}^{(i)}$ having mean $H$ and variance for each $A$, as

$$\frac{1-\alpha}{1+\alpha} H(A)[1 - H(A)].$$

In particular, it follows that the asymptotic distributional behavior of $\nu_{\alpha,\theta+n\alpha}$ is equivalent to that of

$$\sqrt{n} \sum_{i=1}^n \frac{G_{\alpha,i}}{\sum_{j=1}^n G_{\alpha,i}} (P_{\alpha,\alpha}^{(i)} - H)(\cdot).$$

Appealing, again, to the law of large numbers we may instead use the process

$$\frac{1}{\sqrt{n}} \sum_{i=1}^n \frac{G_{\alpha,i}}{\alpha} (P_{\alpha,\alpha}^{(i)} - H)(\cdot)$$

which decomposes into the sum of two asymptotically independent pieces,

$$\frac{1}{\sqrt{n}} \sum_{i=1}^{n} (\frac{G_{\alpha,i}}{\alpha} - 1)(P_{\alpha,\alpha}^{(i)} - H)(\cdot) + \frac{1}{\sqrt{n}} \sum_{i=1}^{n} (P_{\alpha,\alpha}^{(i)} - H)(\cdot).$$

The second term is a measure-like process in the sense of [1] and [39], Section 2.11.1.1. However since we have chosen the class to be VCGC, convergence of this process follows by using arguments similar to that in [41], Section 7.8. The first term then converges by the multiplier CLT in Banach spaces ([24] or [23], Proposition 10.4). $\square$

### 4.1. Dirichlet process asymptotics for $\theta \to \infty$

The next result describes weak convergence of a centered Dirichlet process as $\theta \to \infty$.

**Theorem 4.2.** *Let $\mathcal{F}$ be a VCGC subclass of $L_2(H)$ with envelope $F$ such that $H(F^2) < \infty$. Let, for $\theta > 0$, $P_{0,\theta}$ denote a Dirichlet Process with law $\Pi_{0,\theta}$, having mean $H$, and define $\tau_\theta(\cdot) := \sqrt{\theta}(P_{0,\theta} - H)(\cdot)$. Assume without loss of generality that $\theta = n\kappa$, for $\kappa$ a fixed positive number. Then as $n \to \infty$,*

$$\tau_\theta \rightsquigarrow \mathbb{G}_H \ in \ \ell^\infty(\mathcal{F}),$$

*where $\mathbb{G}_H$ is a $H$-Brownian bridge. Furthermore the limit does not depend on $\kappa$.*

*Proof.* By an argument similar to Proposition 4.1, one can write

$$G_\theta P_{0,\theta} \stackrel{d}{=} \sum_{i=1}^{n} G_{\kappa,i} P_{0,\kappa}^{(i)}.$$

Hence $\tau_\theta$ is asymptotically equivalent to

$$\frac{\sqrt{\kappa}}{\sqrt{n}} \sum_{i=1}^{n} \left( \frac{G_{\kappa,i}}{\kappa} - 1 \right) (P_{0,\kappa}^{(i)} - H)(\cdot) + \frac{\sqrt{\kappa}}{\sqrt{n}} \sum_{i=1}^{n} (P_{0,\kappa}^{(i)} - H)(\cdot).$$

The result then follows analogous to Theorem 4.1. $\square$

### References

[1] ALEXANDER, K. S. (1987). Central limit theorems for stochastic processes under random entropy conditions. *Probab. Theory Related Fields* **75** 351–378. MR0890284

[2] BARRON, A., SCHERVISH, M. J. AND WASSERMAN, L. (1999). The consistency of posterior distributions in nonparametric problems. *Ann. Statist.* **27** 536–561. MR1714718

[3] BICKEL, P. J. AND FREEDMAN, D. A. (1981). Some asymptotic theory for the bootstrap. *Ann. Statist.* **9** 1196–1217. MR0630103

[4] BLACKWELL, D. AND MACQUEEN, J. B. (1973). Ferguson distributions via Pólya urn schemes. *Ann. Statist.* **1** 353–355. MR0362614

[5] BRUNNER, L. J. AND LO, A. Y. (1996). Limiting posterior distributions under mixture of conjugate priors. *Statist. Sinica* **6** 187–197. MR1379056

[6] CHAUMONT, L. and YOR, M. (2003). *Exercises in Probability. A Guided Tour from Measure Theory to Random Processes, via Conditioning.* Cambridge Univ. Press. MR2016344

[7] DAWSON, D. A. and FENG, S. (2006). Asymptotic behavior of the Poisson–Dirichlet distribution for large mutation rate. *Ann. Appl. Probab.* **16** 562–582. MR2244425

[8] FENG, S. (2007). Large deviations for Dirichlet processes and Poisson–Dirichlet distribution with two parameters. *Electron. J. Probab.* **12** 787–807. MR2318410

[9] FERGUSON, T. S. (1973). A Bayesian analysis of some nonparametric problems. *Ann. Statist.* **1** 209–230. MR0350949

[10] FREEDMAN, D. A. and DIACONIS, P. (1983). On inconsistent Bayes estimates in the discrete case. *Ann. Statist.* **11** 1109–1118. MR0720257

[11] GHOSAL, S., GHOSH, J. K. and RAMAMOORTHI, R. V. (1999). Posterior consistency of Dirichlet mixtures in density estimation. *Ann. Statist.* **27** 143–158. MR1701105

[12] GHOSAL, S., GHOSH, J. K. and VAN DER VAART, A. W. (2000). Convergence rates of posterior distributions. *Ann. Statist.* **28** 500–531. MR1790007

[13] GHOSH, J. K. and RAMAMOORTHI, R. V. (2003). *Bayesian Nonparametrics.* Springer, New York. MR1992245

[14] GINÉ, E. (1997). Lectures on some aspects of the bootstrap. *Lectures on Probability Theory and Statistics (Saint-Flour, 1996)* 37–151. *Lecture Notes in Math.* **1665**. Springer, Berlin. MR1490044

[15] GINÉ, E. and ZINN, J. (1984). Some limit theorems for empirical processes (with discussion). *Ann. Probab.* **12** 929–998. MR0757767

[16] GINÉ, E. and ZINN, J. (1990). Bootstrapping general empirical measures. *Ann. Probab.* **18** 851–869. MR1055437

[17] ISHWARAN, H. and JAMES, L. F. (2003). Generalized weighted Chinese restaurant processes for species sampling mixture models. *Statist. Sinica* **13** 1211–1235. MR2026070

[18] ISHWARAN, H. and JAMES, L. F. (2001). Gibbs sampling methods for stick-breaking priors. *J. Amer. Statist. Assoc.* **96** 161–173. MR1952729

[19] JAMES, L. F., LIJOI, A. and PRÜNSTER, I. (2008). Distributions of linear functionals of two-parameter Poisson–Dirichlet random measures. *Ann. Appl. Probab.* **18** 521–551.

[20] JANG, J., LEE, J. and LEE, S. (2007). Posterior consistency of species sampling models. Preprint.

[21] JOYCE, P., KRONE, S. M. and KURTZ, T. G. (2002). Gaussian limits associated with the Poisson–Dirichlet distribution and the Ewens sampling formula. *Ann. Appl. Probab.* **12** 101–124. MR1890058

[22] KIM, Y. and LEE, J. (2004). A Bernstein-von Mises theorem in the nonparametric right-censoring model. *Ann. Statist.* **32** 1492–1512. MR2089131

[23] LEDOUX, M. and TALAGRAND, M. (1991). *Probability in Banach Spaces. Isoperimetry and Processes.* Springer, Berlin. MR1102015

[24] LEDOUX, M. and TALAGRAND, M. (1986). Conditions d'intégrabilité pour les multiplicateurs dans le TLC banachique. *Ann. Probab.* **14** 916–921. MR0841593

[25] LIJOI, A., PRÜNSTER, I. and WALKER, S. G. (2005). On consistency of nonparametric normal mixtures for Bayesian density estimation. *J. Amer. Statist. Assoc.* **100** 1292–1296. MR2236442

[26] LO, A. Y. (1993). A Bayesian bootstrap for censored data. *Ann. Statist.* **21** 100–123. MR1212168

[27] Lo, A. Y. (1987). A large sample study of the Bayesian bootstrap. *Ann. Statist.* **15** 360–375. MR0885742

[28] Lo, A. Y. (1986). A remark on the limiting posterior distribution of the multiparameter Dirichlet process. *Sankhyā Ser. A* **48** 247–249. MR0905464

[29] Lo, A. Y. (1984). On a class of Bayesian nonparametric estimates. I. Density estimates. *Ann. Statist.* **12** 351–357. MR0733519

[30] Lo, A. Y. (1983). Weak convergence for Dirichlet processes. *Sankhyā Ser. A* **45** 105–111. MR0749358

[31] Lo, A. Y. (1982). Bayesian nonparametric statistical inference for Poisson point processes. *Z. Wahrsch. Verw. Gebiete* **59** 55–66. MR0643788

[32] Lynch, J. and Sethuraman, J. (1987). Large deviations for processes with independent increments. *Ann. Probab.* **15** 610–627. MR0885133

[33] Mason, D. M. and Newton, M. A. (1992). A rank statistics approach to the consistency of a general bootstrap. *Ann. Statist.* **20** 1611–1624. MR1186268

[34] Pitman, J. (2006). *Combinatorial Stochastic Processes. Summer School on Probability Theory held in Saint-Flour, July 7–24, 2002. Lecture Notes in Math.* **1875**. Springer, Berlin. MR2245368

[35] Pitman, J. (1996). Some developments of the Blackwell–MacQueen urn scheme. *Statistics, Probability and Game Theory* 245–267. *IMS Lecture Notes Monogr. Ser.* **30**. IMS, Hayward, CA. MR1481784

[36] Pitman, J. and Yor, M. (1997). The two-parameter Poisson–Dirichlet distribution derived from a stable subordinator. *Ann. Probab.* **25** 855–900. MR1434129

[37] Praestgaard, J. and Wellner, J. A. (1993). Exchangeably weighted bootstraps of the general empirical process. *Ann. Probab.* **21** 2053–2086. MR1245301

[38] van der Vaart, A. W. (1998). *Asymptotic Statistics.* Cambridge Univ. Press. MR1652247

[39] van der Vaart, A. W. and Wellner, J. A. (1996). *Weak Convergence and Empirical Processes. With Applications to Statistics.* Springer, New York. MR1385671

[40] Vershik, A., Yor, M. and Tsilevich, N. (2004). On the Markov–Krein identity and quasi-invariance of the gamma process. *J. Math. Sci.* **121** 2303–2310. MR1879060

[41] Ziegler, K. (1997). Functional central limit theorems for triangular arrays of function-indexed processes under uniformly integrable entropy conditions. *J. Multivariate Anal.* **62** 233–272. MR1473875